\documentclass[reqno]{amsart}
\usepackage{amssymb}

\newtheorem{theorem}{Theorem}[section]
\newtheorem{proposition}[theorem]{Proposition}
\newtheorem{lemma}[theorem]{Lemma}
\newtheorem{corollary}[theorem]{Corollary}

\theoremstyle{definition}

\theoremstyle{remark}

\numberwithin{equation}{section}

\begin{document}

\title[Asymptotics of Orthogonal Polynomials]
{Asymptotics of Multivariate Orthogonal Polynomials with
Hyperoctahedral Symmetry}

\author{J.F. van Diejen}
\address{Instituto de Matem\'atica y F\'{\i}sica, Universidad de
Talca, Casilla 747, Talca, Chile}


\thanks{Work supported in part by the Fondo Nacional de Desarrollo
Cient\'{\i}fico y Tecnol\'ogico (FONDECYT) Grant \# 1010217 and by
the Programa Formas Cuadr\'aticas of the Universidad de Talca.}

\begin{abstract}
We present a formula describing the asymptotics of a class of
multivariate orthogonal polynomials with hyperoctahedral symmetry
as the degree tends to infinity. The polynomials under
consideration are characterized by a factorized weight function
satisfying certain analyticity assumptions. As an application, the
large-degree asymptotics of the Koornwinder-Macdonald $BC_N$-type
multivariate Askey-Wilson polynomials is determined.
\end{abstract}

\maketitle

\section{Introduction}\label{sec1}
It is known that the classical families of hypergeometric and
basic hypergeometric orthogonal polynomials form a hierarchy, the
Askey scheme, of which the most general member is given by the
celebrated Askey-Wilson polynomials
\cite{ask-wil:some,koe-swa:askey-scheme}. Other families of
classical orthogonal polynomials, such as e.g. the Hermite,
Laguerre, Jacobi, and Hahn polynomials, all turn out to be special
(limiting) cases of these Askey-Wilson polynomials. Around a
decade ago, Koornwinder introduced a multivariate generalization
of the Askey-Wilson polynomials with hyperoctahedral symmetry
\cite{koo:askey-wilson}, by building upon the pioneering works of
Macdonald on families of orthogonal polynomials associated with
root systems \cite{mac:orthogonal,mac:symmetric,mac:affine}. As it
turns out, these Koornwinder-Macdonald polynomials form again a
master family in the sense that they contain all Macdonald
families associated with the classical root systems as special
cases \cite{koo:askey-wilson,die:commuting}, as well as certain
multivariate versions of the Hermite, Laguerre, Jacobi, and Hahn
polynomials (see e.g. the papers
\cite{bee-opd:certain,bak-for:calogero-sutherland,die:confluent,die:properties}
and references therein). Over the past few years, the properties
of the Koornwinder-Macdonald polynomials have been subject of
investigation in a number of works, leading to a multivariate
generalization of significant part of the theory surrounding the
one-variable Askey-Wilson polynomials
\cite{die:self-dual,die:properties,oko:bc-type,sah:nonsymmetric,die-sto:multivariable,%
nis-kom:algebraic,sto:koornwinder,mim:duality,cha:macdonald,rai:bcn-symmetric}.

A fundamental problem in the theory of orthogonal polynomials is
the question of their asymptotical behavior as the degree tends to
infinity \cite{sze:orthogonal,dei-zho:uniform,dei:orthogonal}. For
the Askey-Wilson polynomials, this asymptotics was determined by
Ismail and Wilson \cite{ism-wil:asymptotic} (leading asymptotics)
and by Ismail \cite{ism:asymptotics} (full asymptotic expansion).
The main purpose of the present work is to lift this asymptotic
analysis to the multivariate level. More specifically, our goal is
to determine the large-degree asymptotics of the
Koornwinder-Macdonald multivariate Askey-Wilson polynomials.  Such
large-degree asymptotics  was determined recently for the
Macdonald polynomials by Ruijsenaars \cite{rui:factorized} (for
the type $A$ root systems) and by the present author
\cite{die:asymptotic} (for arbitrary reduced root systems). The
current work should be seen as an extension of these results to
the case of the Koornwinder-Macdonald polynomials, or from a more
conceptual point of view, as an extension from reduced to
nonreduced root systems. Following the ideas of Ruijsenaars, we in
fact determine the asymptotics of orthogonal polynomials
associated to a fairly large class of weight functions that
factorize in terms of one-dimensional $c$-functions. However,
whereas Ruijsenaars studies homogeneous symmetric polynomials,
here in contrast we consider (Laurent) polynomials in $N$
variables invariant under the action of the hyperoctahedral group
$\Sigma_N\ltimes \mathbb{Z}_2^N$ (thus passing from type $A$ to
type $BC$ root systems). For a specific choice of the
$c$-functions, we end up with the asymptotics of the
Koornwinder-Macdonald polynomials.

It is important to emphasize that---at the multivariate
level---the large-degree asymptotics considered here is not the
only type of asymptotics of interest. Other types of asymptotical
properties of multivariate orthogonal polynomials, involving their
behavior as the number of variables tends to infinity, were for
instance studied by Okounkov and Olshanski for the case of Jack's
hypergeometric degeneration ($q\to 1$) of the Macdonald
polynomials associated with the type $A$ root systems
\cite{oko-ols:asymptotics}.

The paper is organized as follows.  In Section \ref{sec2} we first
define our class of symmetric orthogonal polynomials. An
asymptotic formula for these polynomials is presented in Section
\ref{sec3}.  In Section \ref{sec4} we apply the asymptotic formula
in question to determine the asymptotics of the
Koornwinder-Macdonald polynomials. Finally, Sections \ref{sec5}
and \ref{sec6} wrap up the paper via a series of results
which---when linked together---combine into the proof of the
fundamental asymptotic formula from Section \ref{sec3}.

\section{Multivariate Orthogonal Polynomials}\label{sec2}
\subsection{Symmetric Monomials}
Let $W$ be the hyperoctahedral group given by the semidirect
product of the permutation group $\Sigma_N$ and the N-fold product
of the cyclic group $\mathbb{Z}_2$. The natural action of
$w=(\sigma ,\varepsilon)\in W$ on $\mathbb{R}^N$ is given by
\begin{equation}
\mathbf{x}_w\equiv w(\mathbf{x}) = (\varepsilon_1
x_{\sigma_1},\ldots ,\varepsilon_N x_{\sigma_N})
\end{equation}
(with $\sigma\in \Sigma_N$ and $\varepsilon_j\in \{ 1, -1\}$ for
$j=1,\ldots ,N$). Let $\mathcal{A}^W$ be the algebra of
$W$-invariant trigonometric polynomials on the torus
\begin{equation}
\mathbb{T}_N=\mathbb{R}^N/(2\pi \mathbb{Z})^N.
\end{equation}
The standard basis for $\mathcal{A}^W$ is given by the symmetric
monomials
\begin{subequations}
\begin{equation}
m_\lambda (\mathbf{x}) =\frac{1}{|W_\lambda |} \sum_{w\in W}
e^{i\langle \lambda ,\mathbf{x}_w\rangle},\quad \lambda\in\Lambda,
\end{equation}
where $\langle \lambda
,\mathbf{x}\rangle=\sum_{j=1}^N\lambda_jx_j$,
\begin{equation}
\Lambda = \{ \lambda\in \mathbb{Z}^N \mid
\lambda_1\geq\lambda_2\geq \cdots \geq \lambda_N\geq 0\} ,
\end{equation}
\end{subequations}
and $|W_\lambda |$ denotes the order of the stabilizer subgroup
$W_\lambda=\{ w\in W \mid \lambda_w=\lambda \}$.

\subsection{Orthogonality}
We will partially order the monomial basis $\{ m_\lambda
\}_{\lambda\in \Lambda}$ by means of the hyperoctahedral dominance
order on $\mathbb{Z}^N$:
\begin{equation}\label{po}
\lambda \succeq \mu \Longleftrightarrow \sum_{j=1}^\ell
\lambda_j\geq \sum_{j=1}^\ell \mu_j\quad \text{for}\; \ell
=1,\ldots ,N.
\end{equation}
Let $\Delta (\mathbf{x})$ be an almost everywhere positive weight
function on the torus $\mathbb{T}_N$. We equip $\mathcal{A}^W$
with an inner product structure associated to $\Delta$ via the
definition
\begin{equation}
\langle f , g\rangle_\Delta = \frac{1}{(2\pi
)^N}\int_{\mathbb{T}_N} f(\mathbf{x})
\overline{g(\mathbf{x})}\,\Delta (\mathbf{x}) \text{d} \mathbf{x}
,\qquad \forall f,g\in\mathcal{A}^W
\end{equation}
(where $\overline{g(\mathbf{x})}$ denotes the complex conjugate of
$g(\mathbf{x})$).

Application of the Gram-Schmidt process to the partially ordered
monomial basis produces a basis $\{ P_\lambda \}_{\lambda\in
\Lambda}$ of $\mathcal{A}^W$ of the form
\begin{subequations}
\begin{equation}\label{op1}
P_\lambda (\mathbf{x}) =  \sum_{\mu\in\Lambda,\, \mu \preceq
\lambda} a_{\lambda\mu} m_\mu (\mathbf{x}),\qquad
\lambda\in\Lambda,
\end{equation}
with coefficients $a_{\lambda\mu}\in\mathbb{C}$ such that
\begin{equation}\label{op2}
\langle P_\lambda , m_\mu \rangle_\Delta = 0 \quad\text{if}\; \mu
\prec\lambda \qquad\text{and}\qquad  \langle P_\lambda , P_\lambda
\rangle_\Delta = 1 ,
\end{equation}
\end{subequations}
where we chose $ a_{\lambda \lambda}>0$ by convention. (In other
words, the elements of this new basis are orthogonal when
comparable in the partial order.)

\subsection{Factorized Weight Functions}
From now on we will restrict our attention to a special class of
$W$-invariant weight functions that factorize in terms of
one-dimensional $c$-functions. Specifically, we consider weight
functions of the form
\begin{subequations}
\begin{equation}\label{fw1}
\Delta (\mathbf{x}) =  \frac{1}{|W|\,\mathcal{C}(\mathbf{x})\,
\mathcal{C}(-\mathbf{x})} ,
\end{equation}
with $|W|=2^N N!$ and
\begin{equation}\label{fw2}
\mathcal{C}(\mathbf{x})=\prod_{1\leq j<k\leq N} c_0(x_j+x_k)\,
c_0(x_j-x_k) \prod_{1\leq j\leq N} c_1 (x_j) .
\end{equation}
\end{subequations}
For technical reasons, it will be assumed that the $c$-functions
$c_p(x)$, $p=0,1$ are of the form
\begin{subequations}
\begin{equation}
c_0(x)=(1-e^{-ix})^{-1}\hat{c}_0(e^{-ix})\quad
c_1(x)=(1-e^{-2ix})^{-1}\hat{c}_1(e^{-ix}) .
\end{equation}
Here the reduced $c$-functions $ \hat{c}_{p}(z)$ are taken to be:
analytic on a closed disc $\mathbb{D}_p =\{
z\in\mathbb{C}\mid |z|\leq \varrho_p \}$ of radius $\varrho_p
>1$, zero-free on an open environment of the origin
containing the closed unit disc,
real-valued for $z$ real, and normalized such that $
\hat{c}_{p}(0)=1$. It follows from these conditions that $
\hat{c}_{p}(z)$ and $ 1/ \hat{c}_{p}(z)$ have uniformly
converging Taylor expansions on the closed unit disc of
the form
\begin{equation}\label{taylor-exp}
\hat{c}_p^{\pm 1} (z) = 1 +\sum_{n=1}^\infty a_{n,p}^{\pm}\,
z^n\quad (p=0,1),
\end{equation}
\end{subequations}
with $a^{+}_{n,0}=O(e^{-\epsilon\, n})$ and
$a^{+}_{n,1}=O(e^{-\epsilon\, n/2})$ as $n\to\infty$ for
$$0<\epsilon \leq \min (\log (\varrho_0 ),2\log(\varrho_1)) .$$
Indeed, the asymptotic bound on the Taylor coefficients follows
from the Cauchy formula $a_{n,p}^{+} = \frac{1}{2\pi i}
\oint_{|z|=\varrho_p}\hat{c}_{p}(z) z^{-n-1} \text{d}z$
(whence  $a^{+}_{n,p} =O(e^{-n\log (\varrho_p )})$ as $n$ tends
to $\infty$).

\section{Asymptotic Formulas}\label{sec3}
For $\lambda\in \Lambda$, we define
\begin{equation}\label{asf}
P_\lambda^\infty (\mathbf{x}) = \sum_{w\in W}
\mathcal{C}(\mathbf{x}_w)\, e^{i\langle \lambda
,\mathbf{x}_w\rangle}
\end{equation}
and
\begin{equation}
m(\lambda ) = \min_{j=1,\ldots ,N} \lambda_j-\lambda_{j+1}
\end{equation}
(with the convention that $\lambda_{N+1}\equiv 0$). Let $\| \cdot
\|_\Delta=\sqrt{\langle \cdot ,\cdot \rangle_\Delta}$. The
following theorem states that for $m(\lambda)\to\infty$ the strong
$L^2$-asymptotics  of the polynomials $P_\lambda (\mathbf{x})$ is
given by the functions $P_\lambda^\infty (\mathbf{x})$, with an
exponential error bound governed by the decay rate $\epsilon$ of
the Taylor coefficients $a_{n,p}^+$ of the reduced $c$-functions
$\hat{c}_p(z)$.

\begin{theorem}[Asymptotic Formula A]\label{as1:thm}
One has that
\begin{equation*}
\| P_\lambda -P_\lambda^\infty \|_\Delta = O(e^{-\epsilon\,
m(\lambda )/2})\quad \text{as}\;\; m(\lambda)\to\infty.
\end{equation*}
\end{theorem}

If the polynomials $P_\lambda$ are moreover orthogonal when
non-comparable in the partial order (i.e., if they form an
orthonormal basis of $\mathcal{A}^W$), then we have the following
alternative error bound.

\begin{theorem}[Asymptotic Formula B]\label{as2:thm}
If the basis $\{ P_\lambda\}_{\lambda\in\Lambda}$ is orthogonal,
then one has that
\begin{equation*}
\| P_\lambda -P_\lambda^\infty \|_\Delta = O(\lambda_1^N\,
e^{-\epsilon\, m(\lambda )})\quad \text{as}\;\;
m(\lambda)\to\infty.
\end{equation*}
\end{theorem}

If the growth of $m(\lambda)$ and $\lambda_1$ is proportional,
then the error bound of Theorem \ref{as2:thm} is more efficient
than that of Theorem \ref{as1:thm}. For instance, for $\lambda\in
\Lambda$ fixed and strongly dominant (i.e. with $m(\lambda)>0$),
we get the following asymptotics along the discrete ray
$\lambda\mathbb{N}$.
\begin{corollary}[Ray Asymptotics]\label{ray:cor}
Let $\lambda\in\Lambda$ be fixed and strongly dominant. If the
basis $\{ P_\lambda\}_{\lambda\in\Lambda}$ is orthogonal, then one
has that
\begin{equation*}
\| P_{\ell \lambda} -P_{\ell \lambda}^\infty \|_\Delta = O(\ell^N
e^{-\epsilon\, m(\lambda )\,\ell})\quad \text{as}\;\;
\ell\to\infty.
\end{equation*}
\end{corollary}

In case of polynomial reduced $c$-functions $\hat{c}_p(z)$, the
asymptotic formula turns out to be exact for $m(\lambda)$
sufficiently large.

\begin{theorem}[Exact Asymptotics]\label{exact:thm}
If there exists a nonnegative integer $M$ such that $a_{n,0}^+=0$,
$\forall n > M$ and $a_{n,1}^+=0$, $\forall n > 2M$, then one has
that
\begin{equation*}
P_\lambda (\mathbf{x})={\textstyle
\frac{1}{\mathcal{N}_{\lambda}^{\infty}}} P_\lambda^\infty
(\mathbf{x})
\end{equation*}
for $m(\lambda )\geq M-1$, where $\mathcal{N}_\lambda^\infty\equiv
\| P_\lambda^\infty \|_{\Delta} = 1$ if $m(\lambda ) \geq M $.
\end{theorem}

For $M=0$, the $c$-functions are of the form
$c_0(x)=(1-e^{-ix})^{-1}$ and $c_1(x)=(1-e^{-2ix})^{-1}$,
respectively, and the polynomials $P_\lambda (\mathbf{x})$ amount
in this case to the characters of the symplectic Lie group
$SP(2N;\mathbb{C})$ (with root system $C_N$). The formula of
Theorem \ref{exact:thm} boils then down to the Weyl character
formula.

For $M=1$, the $c$-functions are of the form
$c_0(x)=(1-te^{-ix})(1-e^{-ix})^{-1}$ and $c_1(x)=(1-t_0
e^{-ix})(1-t_1e^{-ix})(1-e^{-2ix})^{-1}$ (with $-1<t,t_0,t_1<1$),
respectively, and the polynomials $P_\lambda (\mathbf{x})$ amount
in this case to Macdonald's generalized Hall-Littlewood
polynomials associated with the root system $BC_N$
\cite{mac:orthogonal}. The formula of Theorem \ref{exact:thm}
boils then down to the standard explicit representation for these
polynomials (cf. Eq. (10.1) of \cite{mac:orthogonal}).

When replacing the partial order $\succeq$ \eqref{po} in the
Gram-Schmidt process of definition \eqref{op1}, \eqref{op2} by the
lexicographical ordering $\dot\succeq$, one ends up with an
orthonormal basis $\{ \dot{P}_\lambda \}_{\lambda\in\Lambda}$ of
$\mathcal{A}^W$. It is clear from the proofs in Sections
\ref{sec5} and \ref{sec6} that the asymptotics of the polynomials
$\dot{P}_\lambda (\mathbf{x})$, $\lambda\in\Lambda$ is again given
by Theorem \ref{as1:thm}, by Theorem \ref{as2:thm} (and thus
Corollary \ref{ray:cor}), and by Theorem \ref{exact:thm} (with the
same asymptotic functions $P^\infty_\lambda (\mathbf{x})$ and the
same error bound $\epsilon$). (The crux is that Proposition
\ref{orthogonality:prp} and Lemma \ref{dimbound:lem} below remain
valid when replacing the partial order $\succeq$ by the
lexicographical linear refinement $\dot\succeq$.) In particular,
it follows from this observation that in the situation of Theorem
\ref{exact:thm}, one has that $\langle P_\lambda
,P_\mu\rangle_\Delta =\delta_{\lambda\mu}$ for
$\lambda,\mu\in\Lambda$ such that $\min (m(\lambda),m(\mu))\geq
M-1$ (even if $\lambda$ and $\mu$ are not comparable in the
partial order $\succeq$ \eqref{po}).

\section{Specialization to Koornwinder-Macdonald
Polynomials}\label{sec4} By picking reduced $c$-functions
$\hat{c}_p(z)$, $p=0,1$, of the form
\begin{subequations}
\begin{equation}\label{krc-f1}
\hat{c}_0(z)=\frac{(tz;q)_\infty}{(qz;q)_\infty},\qquad
\hat{c}_1(z)=\frac{\prod_{r=0}^3(t_rz;q)_\infty}{(qz^2;q)_\infty},
\end{equation}
where $(z;q)_\infty\equiv\prod_{n=0}^\infty (1-zq^n)$, and with
parameters subject to the constraints
\begin{equation} \label{krc-f2}
0<q<1,\qquad
-1<t,t_r< 1\quad (r=0,\ldots ,3),
\end{equation}
\end{subequations}
the weight
function $\Delta (\mathbf{x})$ \eqref{fw1}, \eqref{fw2}
specializes to
\begin{subequations}
\begin{equation}\label{kw1}
\Delta (\mathbf{x}) = \frac{1}{2^N
N!\,\mathcal{C}(\mathbf{x})\,\mathcal{C}(-\mathbf{x})} ,
\end{equation}
with
\begin{eqnarray}
\mathcal{C} (\mathbf{x})&= &\prod_{1\leq j < k\leq N}
\frac{(t\,e^{-i(x_j+x_k)},t\,e^{-i(x_j-x_k)};q)_\infty}
{(e^{-i(x_j+x_k)},e^{-i(x_j-x_k)};q)_\infty} \nonumber \\
&& \times \prod_{1\leq j\leq N} \frac{\prod_{r=0}^3(t_r\,
e^{-ix_j};q)_\infty}{(e^{-2ix_j} ;q)_\infty}  \label{kw2}
\end{eqnarray}
\end{subequations}
(and $(z_1,z_2,\ldots ,z_k;q)_\infty\equiv (z_1;q)_\infty
(z_2,q)_\infty\cdots (z_k;q)_\infty$). It was shown by Koornwinder
that the polynomials $P_\lambda (\mathbf{x})$, $\lambda\in\Lambda$
\eqref{op1}, \eqref{op2}, associated to the weight function
$\Delta (\mathbf{x})$ \eqref{kw1}, \eqref{kw2}, form an orthogonal
system \cite{koo:askey-wilson} (i.e., the polynomials are also
orthogonal when non-comparable in the partial order). The
conditions on the parameters $q$, $t$ and $t_0,\ldots ,t_3$ in
Equation \eqref{krc-f2} ensure that the reduced $c$-functions
$\hat{c}_0(z)$, $\hat{c}_1(z)$ in Equation \eqref{krc-f1} satisfy
the technical requirements stipulated in Section \ref{sec2}. In
particular, for the bound on the decay rate of the Taylor
coefficients $a_{n,p}^+$ we have $
\varrho_0=q^{-1}$ and  $\varrho_1=q^{-1/2}$, so we may choose
(any) $\epsilon \in (0, \log (1/q))$.

Specialization of the results of Section \ref{sec3} to the weight
function $\Delta (\mathbf{x})$ \eqref{kw1}, \eqref{kw2}
immediately entails the main application of our asymptotic
analysis.

\begin{corollary}[Asymptotics of Koornwinder-Macdonald
Polynomials]\label{kas:cor} The asymptotics of the
Koornwinder-Macdonald polynomials is governed by Theorem
\ref{as1:thm}, Theorem \ref{as2:thm}, and Corollary \ref{ray:cor},
with asymptotic functions $P_\lambda^\infty(\mathbf{x})$
\eqref{asf} characterized by the product $c$-function
$\mathcal{C}(\mathbf{x})$ \eqref{kw2}, and an error bound with a
decay rate that is at least as fast as any $\epsilon$ taken from
the interval $(0,\log (1/q))$.
\end{corollary}

\section{The Asymptotic Functions}\label{sec5}
This section exhibits some properties of the asymptotic functions
$P_\lambda^\infty (\mathbf{x})$ \eqref{asf} that are needed in the
proof of the asymptotic formulas stated in Section \ref{sec3}.

\begin{proposition}[Partial Biorthogonality]\label{orthogonality:prp}
Let $\lambda ,\mu\in\Lambda$ with $\mu\preceq\lambda$. Then
\begin{equation*}
\langle P_\lambda^\infty, m_\mu  \rangle_\Delta =
\begin{cases}
0 & \text{if}\; \mu\prec\lambda ,\\
1 & \text{if}\; \mu =\lambda .
\end{cases}
\end{equation*}
\end{proposition}
\begin{proof}
From the $W$-invariance of the weight function $\Delta
(\mathbf{x})=\frac{1}{\mathcal{C}(\mathbf{x})\mathcal{C}(-\mathbf{x})}$
it is clear that
\begin{eqnarray*}
\lefteqn{\langle P_\lambda^\infty, m_\mu  \rangle_\Delta } && \\
&=&\!\! \frac{1 }{ (2\pi)^N |W|\, |W_\mu |} \int_{\mathbb{T}_N}
\frac{1}{\mathcal{C}(\mathbf{x})\mathcal{C}(-\mathbf{x})}
\sum_{w_1\in W} \mathcal{C}(\mathbf{x}_{w_1}) e^{i\langle \lambda
,\mathbf{x}_{w_1}\rangle}\sum_{w_2\in W}e^{-i\langle \mu
,\mathbf{x}_{w_2}\rangle} \text{d}\mathbf{x} \\
&=&\!\! \frac{1 }{ (2\pi)^N\, |W_\mu |} \sum_{w\in W}
\int_{\mathbb{T}_N} \frac{1}{\mathcal{C}(-\mathbf{x})} e^{i\langle
\lambda - \mu_w ,\mathbf{x}\rangle} \text{d}\mathbf{x}.
\end{eqnarray*}
The integral on the last line picks up the constant term of the
integrand (times $(2\pi )^N$). It is immediate from our
assumptions on the structure of the $c$-functions that
$1/\mathcal{C}(-\mathbf{x})$ has a (uniformly converging) Fourier
expansion of the form $1+\sum_{\mathbf{n}\in\mathbb{Z}^N,\,
\mathbf{n}\succ \mathbf{0}} \mathbf{a}_{\mathbf{n}}^- e^{i\langle
\mathbf{n} ,\mathbf{x}\rangle}$, whence the constant term in
question is equal to $1$ if $\lambda =\mu_w=\mu$ and equal to $0$
otherwise. (Here we used the standard fact that $\mu\succeq \mu_w$
for all $\mu \in\Lambda$ and $w\in W$, cf. also Lemma
\ref{saturation:lem} below.)
\end{proof}

By factoring-off the denominators of the $c$-functions, one
rewrites $P_\lambda^\infty (\mathbf{x})$ \eqref{asf} as
\begin{subequations}\label{asf2}
\begin{equation}
P_\lambda^\infty (\mathbf{x}) = \delta^{-1}(\mathbf{x}) \sum_{w\in
W} \det (w) \hat{\mathcal{C}}(\mathbf{x}_w) e^{i\langle \lambda
+\rho , \mathbf{x}_w \rangle},
\end{equation}
with
\begin{eqnarray}
\hat{\mathcal{C}}(\mathbf{x})&=&\prod_{1\leq j<k\leq N}
\hat{c}_0(e^{-i(x_j+x_k)})\, \hat{c}_0(e^{-i(x_j-x_k)})
\prod_{1\leq j\leq N} \hat{c}_1 (e^{-ix_j}) , \\ \delta
(\mathbf{x}) &=&\prod_{1\leq j<k\leq N} \bigl(
e^{i(x_j+x_k)/2}-e^{-i(x_j+x_k)/2}   \bigr)
(e^{i(x_j-x_k)/2}-e^{-i(x_j-x_k)/2}) \nonumber \\
&& \times \prod_{1\leq j\leq N} (e^{ix_j}-e^{-ix_j}) ,
\end{eqnarray}
and
\begin{equation}
\rho= \sum_{j=1}^N (N+1-j) \mathbf{e}_j.
\end{equation}
\end{subequations}
(Here $\mathbf{e}_j$ denotes the $j^{th}$ unit vector in the
standard basis of $\mathbb{R}^N$.) We introduce the following
polynomial truncation of the asymptotic function $P_\lambda^\infty
(\mathbf{x})$:
\begin{subequations}
\begin{equation}\label{asfm}
P_\lambda^{(m)} (\mathbf{x}) = \delta^{-1}(\mathbf{x}) \sum_{w\in
W} \det (w) \hat{\mathcal{C}}^{(m)}(\mathbf{x}_w) e^{i\langle
\lambda +\rho , \mathbf{x}_w \rangle},
\end{equation}
with
\begin{equation}
\hat{\mathcal{C}}^{(m)}(\mathbf{x})=\prod_{1\leq j<k\leq N}
\hat{c}_0^{(m)}(e^{-i(x_j+x_k)})\,
\hat{c}_0^{(m)}(e^{-i(x_j-x_k)}) \prod_{1\leq j\leq N}
\hat{c}_1^{(m)} (e^{-ix_j}) ,
\end{equation}
where $\hat{c}_0^{(m)} (z)$ and $\hat{c}_1^{(m)} (z)$ consists of
the first $m+1$ and $2m+1$ terms of the Taylor expansions of the
reduced $c$-functions $\hat{c}_0 (z)$ and $\hat{c}_1 (z)$,
respectively, i.e.
\begin{equation}
\hat{c}_0^{(m)} (z) =1 +\sum_{n=1}^m a_{n,0}^{+}\, z^n,\qquad
\hat{c}_1^{(m)} (z) =1 +\sum_{n=1}^{2m} a_{n,1}^{+}\, z^n
\end{equation}
\end{subequations}
(with the coefficients $a_{n,p}^+$ ($p=0,1$) being defined by Eq.
\eqref{taylor-exp}).

\begin{proposition}[Asymptotic Error Bound]\label{as-bound:prp}
One has that
\begin{equation*}
P_\lambda^\infty (\mathbf{x}) =
P_\lambda^{(m)}(\mathbf{x})+\mathcal{E}_\lambda^{(m)}(\mathbf{x}),
\end{equation*}
with
\begin{equation*}
\| \mathcal{E}_\lambda^{(m)} \|_\Delta = O(e^{-\epsilon\, m})
\qquad \text{as}\;\; m \longrightarrow\infty
\end{equation*}
(uniformly in $\lambda$).
\end{proposition}
\begin{proof}
Let us write
\begin{equation*}
\hat{\mathcal{C}}(\mathbf{x}) =\hat{\mathcal{C}}^{(m)}(\mathbf{x})
+ \mathcal{R}^{(m)}(\mathbf{x}).
\end{equation*}
The error between $P_\lambda^{(m)}(\mathbf{x})$ and
$P_\lambda^\infty (\mathbf{x})$ is then given by (cf. Eqs.
\eqref{asf2}, \eqref{asfm})
\begin{equation*}
\mathcal{E}_\lambda^{(m)} (\mathbf{x}) = \delta^{-1}(\mathbf{x})
\sum_{w\in W} \det (w) \mathcal{R}^{(m)}(\mathbf{x}_w) e^{i\langle
\lambda +\rho , \mathbf{x}_w \rangle} .
\end{equation*}
The quotient $\Delta (\mathbf{x})/|\delta (\mathbf{x})|^2=
1/(\hat{\mathcal{C}}(\mathbf{x})\hat{\mathcal{C}}(-\mathbf{x}))$
is smooth on the torus $\mathbb{T}_N$ due to the absence of zeros
in the $c$-functions. Hence, to prove the error bound on $\|
\mathcal{E}_\lambda^{(m)} \|_\Delta$ it is enough to show that
$\max_{\mathbf{x}\in\mathbb{T}_N} (\mathcal{R}^{(m)}(\mathbf{x}))=
O(e^{-\epsilon\, m})$. To this end we set $\hat{c}_{p}(z) =
\hat{c}_{p}^{(m)}(z)+r_{p}^{(m)}(z)$, whence
\begin{eqnarray*}
\mathcal{R}^{(m)}(\mathbf{x})&=& \prod_{1\leq j<k\leq N}
\left(\hat{c}_0^{(m)}(e^{-i(x_j+x_k)})+r_0^{(m)}(e^{-i(x_j+x_k)})\right)
\\
&& \makebox[3em]{} \times\left(
\hat{c}_0^{(m)}(e^{-i(x_j-x_k)})+r_0^{(m)}(e^{-i(x_j-x_k)})\right)
\\&& \times\prod_{1\leq j\leq N} \left( \hat{c}_1^{(m)}
(e^{-ix_j})+r_1^{(m)} (e^{-ix_j})\right) \\
&-& \prod_{1\leq j<k\leq N} \hat{c}_0^{(m)}(e^{-i(x_j+x_k)})\,
\hat{c}_0^{(m)}(e^{-i(x_j-x_k)}) \prod_{1\leq j\leq N}
\hat{c}_1^{(m)} (e^{-ix_j}) .
\end{eqnarray*}
The bound on $\max_{\mathbf{x}\in\mathbb{T}_R}
(\mathcal{R}^{(m)}(\mathbf{x}))$ thus follows since
$\max_{|z|=1}(\hat{c}^{(m)}_{p}(z))=O(1)$ and
$\max_{|z|=1}(r^{(m)}_{p}(z))=O(e^{-\epsilon\, m})$, in view of
the $O(e^{-\epsilon\, n})$ and $O(e^{-\epsilon\, n/2})$ decay
rates of the expansion coefficients $a_{n,0}^+$ and $a_{n,1}^+$
for the reduced $c$-functions $\hat{c}_{0}(z)$ and
$\hat{c}_{1}(z)$, respectively.
\end{proof}

For $m$ sufficiently small, the polynomial truncation
$P_\lambda^{(m)} (\mathbf{x})$ of the asymptotic function
$P_\lambda^{\infty} (\mathbf{x})$ expands triangularly on the
monomial basis. This observation hinges on the following lemma.

\begin{lemma}\label{saturation:lem}
Let $\lambda\in\Lambda$ and let
\begin{equation*}
\mu=\lambda -\sum_{1\leq j<k\leq N} \left(
n_{jk}^+(\mathbf{e}_j+\mathbf{e}_k)+n_{jk}^-(\mathbf{e}_j-\mathbf{e}_k)\right)+
\sum_{1\leq j\leq N} n_j\mathbf{e}_j ,
\end{equation*}
with  $0\leq n_{jl}^+,n_{jk}^-\leq m(\lambda)$ and $0\leq n_j\leq
2 m(\lambda)$. Then one has that
\begin{equation*}
\mu_w\preceq  \lambda,\quad \forall w\in W.
\end{equation*}
Furthermore, if  $0\leq n_{jl}^+,n_{jk}^-< m(\lambda)$ and $0\leq
n_j< 2 m(\lambda)$, then the equality $\mu_w=\lambda$ is assumed
if and only if $w=\text{Id}$ (and all $n_{jk}^+,n_{jk}^-,n_j$
vanish).
\end{lemma}
\begin{proof}
The components of $\mu$ are given by
\begin{equation*}
\mu_j=\lambda_j-n_j-\sum_{1\leq k <j} (n_{kj}^+
-n_{kj}^-)-\sum_{j< k <N} (n_{jk}^+ +n_{jk}^-),
\end{equation*}
$j=1,\ldots ,N$. It thus follows that
\begin{eqnarray}
\sum_{j=1}^\ell \varepsilon_j\, \mu_{\sigma_j}\!\! &= & \!\!
\sum_{j\in J_+} \Bigl( \lambda_{\sigma_j}-n_{\sigma_j}-\sum_{1\leq
k <\sigma_j} (n_{k\sigma_j}^+ -n_{k\sigma_j}^-)-\sum_{\sigma_j< k
<N} (n_{\sigma_jk}^+
+n_{\sigma_jk}^-)\Bigr)  +\nonumber\\
&& \!\! \sum_{j\in J_-} \Bigl(
n_{\sigma_j}-\lambda_{\sigma_j}+\sum_{1\leq k <\sigma_j}
(n_{k\sigma_j}^+
-n_{k\sigma_j}^-)+\sum_{\sigma_j< k <N} (n_{\sigma_jk}^+ +n_{\sigma_jk}^-)\Bigr)\nonumber \\
&\leq & \!\! \sum_{j\in J_+}\Bigl( \lambda_{\sigma_j} +
\sum_{\begin{subarray}{c} 1\leq k< \sigma_j \\ k\not\in \sigma
(J_+)
\end{subarray}} n_{k\sigma_j}^- \Bigr) + \nonumber\\
&& \!\! \sum_{j\in J_-}  \Bigl( n_{\sigma_j}-\lambda_{\sigma_j} +
\sum_{\begin{subarray}{c} 1\leq k< \sigma_j \\ k\not\in \sigma
(J_+)
\end{subarray}} n_{k\sigma_j}^+ +
\sum_{\begin{subarray}{c} \sigma_j<k\leq N \\ k\not\in \sigma
(J_+)
\end{subarray}} n_{\sigma_jk}^+ +
\sum_{\begin{subarray}{c} \sigma_j <k\leq N\\ k\not\in \sigma
(J_-)
\end{subarray}} n_{\sigma_jk}^-
\Bigr) ,    \label{i1}
\end{eqnarray}
where  $J_+=\{ 1\leq j\leq \ell \mid \varepsilon_j= +1\}$ and
$J_-=\{ 1\leq j\leq \ell \mid \varepsilon_j= -1\}$. The proof of
the first part of the lemma now hinges on successive application
of the following three elementary `transportation' inequalities
\begin{equation*}\label{transport}
\lambda_j+n\leq \lambda_{j-1}\quad \text{(A)},\qquad
n-\lambda_j\leq -\lambda_{j+1}\quad\text{(B)},\qquad m
-\lambda_N\leq \lambda_N\quad\text{(C)},
\end{equation*}
for $0\leq n\leq m(\lambda)$ and $0\leq m\leq 2m(\lambda)$.
Indeed, ordering of the components $\lambda_{\sigma_j}$, $j\in
J_+$ from small to large and iterated application of inequality
(A) (so as to `transport' to $\lambda_1,\lambda_2,\ldots
,\lambda_{|J_+|}$, respectively) readily entails that
\begin{equation}\label{i2}
\sum_{j\in J_+}\Bigl( \lambda_{\sigma_j} +
\sum_{\begin{subarray}{c} 1\leq k< \sigma_j \\ k\not\in \sigma
(J_+)
\end{subarray}} n_{k\sigma_j}^- \Bigr) \leq \lambda_1
+\lambda_2+\cdots +\lambda_{|J_+|}
\end{equation}
(where $|J_+|$ denotes the number of elements of $J_+$). In a
similar way we obtain that
\begin{eqnarray}
\lefteqn{\sum_{j\in J_-}  \Bigl( n_{\sigma_j}-\lambda_{\sigma_j} +
\sum_{\begin{subarray}{c} 1\leq k< \sigma_j \\ k\not\in \sigma
(J_+)
\end{subarray}} n_{k\sigma_j}^+ +
\sum_{\begin{subarray}{c} \sigma_j<k\leq N \\ k\not\in \sigma
(J_+)
\end{subarray}} n_{\sigma_jk}^+ +
\sum_{\begin{subarray}{c} \sigma_j <k\leq N\\ k\not\in \sigma
(J_-)
\end{subarray}} n_{\sigma_jk}^-
\Bigr) }&& \nonumber \\
&\stackrel{(i)}\leq & \sum_{j\in J_-}  \Bigl( n_{\sigma_j} +
\sum_{\begin{subarray}{c} 1\leq k< \sigma_j \\ k\not\in \sigma
(J_+)
\end{subarray}} n_{k\sigma_j}^+ +
\sum_{\begin{subarray}{c} \sigma_j<k\leq N \\ k\not\in \sigma
(J_+)
\end{subarray}} n_{\sigma_jk}^+
\Bigr)  \nonumber \\
&& -(\lambda_{N-|J_-|+1} +\lambda_{N-|J_-|+2}+\cdots
+\lambda_{N})  \nonumber \\
&\stackrel{(ii)}{\leq}& \lambda_{|J_+|+1}+\lambda_{|J_+|+2}+\cdots
+\lambda_{|J_+|+|J_-|} .   \label{i3}
\end{eqnarray}
Here the inequality {\em (i)} is inferred by ordering
$-\lambda_{\sigma_j}$, $j\in J_-$ from large to small, followed by
iterated application of inequality (B) (`transporting' to
$-\lambda_{N+1-j}$, $j=1,\ldots ,|J_-|$, respectively); the
inequality {\em (ii)} then follows by iterated application of,
respectively, inequality (B) (`transporting' from
$-\lambda_{N+1-j}$ to $-\lambda_N$), inequality (C) (`flipping'
the sign from $-\lambda_N$ to $+\lambda_N$), and inequality (A)
(`transporting' back to $\lambda_{|J_+| +j}$) (with $j=1,\ldots ,
|J_-|$). Combination of Eqs. \eqref{i1}, \eqref{i2} and \eqref{i3}
now gives that $\sum_{j=1}^\ell \varepsilon_j\, \mu_{\sigma_j}\leq
\sum_{j=1}^\ell \lambda_j$ for $\ell=1,\ldots ,N$, which proves
the first part of the lemma. To prove the second part, one
observes that the elementary inequalities (A), (B), and (C),
become strict for $0\leq n< m(\lambda )$ and $0\leq m<
2m(\lambda)$. Hence, for $0\leq n_{jl}^+,n_{jk}^-< m(\lambda)$ and
$0\leq n_j< 2 m(\lambda)$ the inequality in Eq. \eqref{i2} becomes
strict unless $\sigma (J_+)=\{ 1,2, \ldots , |J_+|\}$, and the
inequality in Eq. \eqref{i3} becomes strict unless $J_-=\emptyset$
(so $J_+=\{ 1,2,\ldots ,\ell\}$). Thus, the upshot is that now
$\sum_{j=1}^\ell \varepsilon_j\, \mu_{\sigma_j}= \sum_{j=1}^\ell
\lambda_j$ for $\ell =1,\ldots ,N$ if and only if $\{
\sigma_1,\ldots ,\sigma_\ell\} = \{ 1,\ldots ,\ell\}$ for $\ell =
1,\ldots ,N$ and $\varepsilon_j=+1$ for $j=1,\ldots ,N$, i.e., if
and only if $w=\text{Id}$ (whence all $n_{jl}^+,n_{jk}^-$ and
$n_j$ vanish).
\end{proof}

\begin{proposition}[Triangularity]\label{triangular:prp}
For $\lambda\in\Lambda$ and $m\leq m(\lambda)+1$, one has that
\begin{equation*}
P_\lambda^{(m)} (\mathbf{x}) = \sum_{\mu\in\Lambda,\, \mu \preceq
\lambda} a_{\lambda\mu}^{(m)} m_\mu (\mathbf{x}),
\end{equation*}
with $a_{\lambda\mu}^{(m)}\in\mathbb{C}$. Furthermore, for $m\leq
m(\lambda)$ the polynomial $P_\lambda^{(m)} (\mathbf{x})$ is monic
(i.e. $a_{\lambda\lambda}^{(m)}=1$).
\end{proposition}
\begin{proof}
It is immediate from the definition that the truncated asymptotic
function $P_\lambda^{(m)} (\mathbf{x})$ \eqref{asfm} can be
written as a finite linear combination of symmetric polynomials of
the form
\begin{subequations}
\begin{equation}\label{terms}
\delta^{-1}(\mathbf{x}) \sum_{w\in W} \det (w)  e^{i\langle
\lambda +\rho -\mathbf{n}, \mathbf{x}_w \rangle} ,
\end{equation}
where
\begin{equation}
\mathbf{n}=\sum_{1\leq j<k\leq N} \left(
n_{jk}^+(\mathbf{e}_j+\mathbf{e}_k)+n_{jk}^-(\mathbf{e}_j-\mathbf{e}_k)\right)+
\sum_{1\leq j\leq N} n_j\mathbf{e}_j,
\end{equation}
\end{subequations}
with $0\leq n_{jl}^+,n_{jk}^-\leq m$ and $0\leq n_j\leq 2 m$. The
polynomial in Eq. \eqref{terms} vanishes if $\lambda +\rho
-\mathbf{n}$ is singular with respect to the action of the Weyl
group (i.e. if it has a nontrivial stabilizer) and otherwise it is
equal, possibly up to a sign, to the Weyl character
\begin{equation}
\chi_\mu =\delta^{-1}(\mathbf{x}) \sum_{w\in W} \det (w)
e^{i\langle \mu +\rho , \mathbf{x}_w \rangle},
\end{equation} where $\mu$ is the unique dominant weight
in the translated Weyl orbit $W(\lambda +\rho -\mathbf{n})-\rho$.
It follows from the first part of Lemma \ref{saturation:lem}, and
the assumption $m\leq m(\lambda)+1=m(\lambda +\rho)$, that $\mu
\preceq \lambda$. It moreover follows from the second part of
Lemma \ref{saturation:lem} that for $m\leq m(\lambda)$, one has
that $\mu=\lambda$ if and only if all $n_{jl}^+$, $n_{jk}^-$ and
$n_j$ are zero. The proposition now follows from the well-known
fact that the Weyl characters are monic $W$-invariant polynomials
that expand triangularly on the basis of monomial symmetric
functions.
\end{proof}

We conclude this section with estimates for the norm of the
asymptotic function $P_\lambda^\infty (\mathbf{x})$ and the for
the leading coefficient in the monomial expansion of the
normalized polynomial $P_\lambda (\mathbf{x})$.

\begin{proposition}[Norm Estimate]\label{normbound1:prp}
One has that
\begin{equation*}
\| P_\lambda^\infty \|_\Delta = 1 +O(e^{-\epsilon\, m (\lambda)})
\qquad \text{as}\;\; m(\lambda)\longrightarrow\infty.
\end{equation*}
\end{proposition}
\begin{proof}
It is clear that
\begin{equation*}
\| P_\lambda^\infty \|_\Delta^2 =\langle
P_\lambda^\infty,P_\lambda^\infty \rangle_\Delta = \langle
P_\lambda^\infty,P_\lambda^{(m(\lambda ))}
 \rangle_\Delta +\langle P_\lambda^\infty,\mathcal{E}_\lambda^{(m(\lambda ))}
 \rangle_\Delta .
\end{equation*}
The proposition now follows from the observation that $\langle
P_\lambda^\infty,P_\lambda^{(m(\lambda ))}
 \rangle_\Delta=1$ by Propositions \ref{orthogonality:prp} and
 \ref{triangular:prp},
combined with the error estimate
\begin{equation*} | \langle
P_\lambda^\infty,\mathcal{E}_\lambda^{(m(\lambda ))}
 \rangle_\Delta |\leq \| P_\lambda^\infty \|_\Delta \| \mathcal{E}_\lambda^{(m(\lambda ))}  \|_\Delta
=O(e^{-\epsilon\, m (\lambda)}),
\end{equation*}
by Proposition \ref{as-bound:prp} (using also that
$\max_{\mathbf{x}\in\mathbb{T}_N} | \delta(\mathbf{x})
P_\lambda^\infty (\mathbf{x}) |\leq |W|
\max_{\mathbf{x}\in\mathbb{T}_N} | \hat{\mathcal{C}} (\mathbf{x})
|$ and that $\Delta (\mathbf{x})/|\delta (\mathbf{x})|^2$ is
bounded on $\mathbb{T}_N$).
\end{proof}

Let $\mathcal{N}_\lambda =a_{\lambda\lambda}^{-1}$, where
$a_{\lambda\lambda}$ ($>0$) represents the leading coefficient of
the polynomial $P_\lambda (\mathbf{x})$ \eqref{op1}, \eqref{op2}
in the monomial basis.

\begin{proposition}[Leading Coefficient]\label{normbound2:prp}
One has that
\begin{equation*}
\mathcal{N}_{\lambda}=1+O(e^{-\epsilon\, m (\lambda)}) \qquad
\text{as}\;\; m(\lambda)\longrightarrow\infty.
\end{equation*}

\end{proposition}
\begin{proof}
A sequence of elementary manipulations entails that
\begin{eqnarray*}
a_{\lambda\lambda} &\stackrel{(i)}{=}& \langle P_\lambda ,
P_\lambda^\infty\rangle_\Delta \\
&\stackrel{(ii)}{=}& \langle P_\lambda ,
P_\lambda^{(m(\lambda))}+\mathcal{E}_\lambda^{(m(\lambda))} \rangle_\Delta \\
&\stackrel{(iii)}{=}& \langle P_\lambda ,
a_{\lambda\lambda}^{-1}P_\lambda+\mathcal{E}_\lambda^{(m(\lambda))} \rangle_\Delta \\
&\stackrel{(iv)}{=}& a_{\lambda\lambda}^{-1}+\langle P_\lambda ,
\mathcal{E}_\lambda^{(m(\lambda))}
\rangle_\Delta \\
&\stackrel{(v)}{=}& a_{\lambda\lambda}^{-1}+O(e^{-\epsilon\, m
(\lambda)}),
\end{eqnarray*}
whence
$\mathcal{N}_\lambda^{-1}=a_{\lambda\lambda}=1+O(e^{-\epsilon\, m
(\lambda)})$. Here we used respectively {\em (i)} Eq.~\eqref{op1}
and Proposition \ref{orthogonality:prp}, {\em (ii)} Proposition
\ref{as-bound:prp}, {\em (iii)} Proposition \ref{triangular:prp}
and Eqs. \eqref{op1}, \eqref{op2}, {\em (iv)} Eq.~\eqref{op2}, and
{\em (v)} the Cauchy-Schwarz inequality and Proposition
\ref{as-bound:prp}.
\end{proof}

\section{Proofs of the Main Theorems}\label{sec6}
By combining the properties in Section \ref{sec5}, we arrive at
the proofs of the theorems stated in Section \ref{sec3}.

\subsection{Proof of Theorem \ref{as1:thm}}
Straightforward manipulations reveal that
\begin{eqnarray*}
\| P_\lambda -P_\lambda^\infty\|_\Delta^2 &=&  \langle P_\lambda,
P_\lambda \rangle_\Delta - \langle P_\lambda, P_\lambda^\infty
\rangle_\Delta- \langle P_\lambda^\infty, P_\lambda \rangle_\Delta
+ \langle
P_\lambda^\infty, P_\lambda^\infty \rangle_\Delta \\
&\stackrel{(i)}{=}& 1 -2\, \mathcal{N}_\lambda^{-1}+\| P_\lambda^\infty\|_\Delta^2 \\
&\stackrel{(ii)}{=}&  O(e^{-\epsilon\, m (\lambda)}),
\end{eqnarray*}
whence $\| P_\lambda -P_\lambda^\infty\|_\Delta = O(e^{-\epsilon\,
m (\lambda)/2})$. Step {\em (i)} hinges on Eq. \eqref{op1} and
Proposition \ref{orthogonality:prp}, which implies that $\langle
P_\lambda, P_\lambda^\infty \rangle_\Delta=\langle
P_\lambda^\infty, P_\lambda
\rangle_\Delta=a_{\lambda\lambda}=\mathcal{N}_\lambda^{-1}$. Step
{\em (ii)} follows from the estimates in Proposition
\ref{normbound1:prp} and Proposition \ref{normbound2:prp}.

\subsection{Proof of Theorem \ref{as2:thm}}
From Proposition \ref{triangular:prp} and Eqs. \eqref{op1},
\eqref{op2} it is immediate that for $m\leq m(\lambda)$
\begin{equation}\label{Pm-orthogonalexp}
P_\lambda^{(m)} (\mathbf{x}) = \mathcal{N}_\lambda P_\lambda
(\mathbf{x})+ \sum_{\mu\in\Lambda,\, \mu \prec \lambda}
b_{\lambda\mu}^{(m)} P_\mu (\mathbf{x}),
\end{equation}
where $\mathcal{N}_\lambda=a_{\lambda\lambda}^{-1}$ and
$b_{\lambda\mu}^{(m)}\in\mathbb{C}$.

\begin{lemma}\label{termbound:lem}  If the basis $\{ P_\lambda
\}_{\lambda\in\Lambda}$ is {\em orthonormal}, then
\begin{equation*}
| b_{\lambda\mu}^{(m)}| \leq \| \mathcal{E}_\lambda^{(m)}
\|_\Delta
\end{equation*}
for $\lambda,\mu\in\Lambda$ with $\mu \prec\lambda$ and $m\leq
m(\lambda)$.
\end{lemma}
\begin{proof}
Respectively applying Eq. \eqref{Pm-orthogonalexp}, Proposition
\ref{as-bound:prp}, and Proposition \ref{orthogonality:prp},
readily entails that
\begin{equation*}
b_{\lambda\mu}^{(m)} = \langle P_\lambda^{(m)},P_\mu\rangle_\Delta
= \langle P_\lambda^{\infty}-
\mathcal{E}_\lambda^{(m)},P_\mu\rangle_\Delta
 = -\langle
\mathcal{E}_\lambda^{(m)},P_\mu\rangle_\Delta .
\end{equation*}
Hence $| b_{\lambda\mu}^{(m)}| \leq \| \mathcal{E}_\lambda^{(m)}
\|_\Delta$ by the Cauchy-Schwarz inequality.
\end{proof}

\begin{lemma}\label{dimbound:lem} For $\lambda\in\Lambda$ let
$\mathcal{A}_\lambda^W=\text{Span}\{ m_\mu \}_{\mu\in\Lambda,\,
\mu\preceq\lambda}$. Then
\begin{equation*}
\dim (\mathcal{A}_\lambda^W) \leq \bigl( 1+\lambda_1 \bigr)^N .
\end{equation*}
\end{lemma}
\begin{proof}
Immediate from the observation that for $\mu\in\Lambda$ the
inequality $\mu\preceq\lambda$ implies that $0\leq \mu_j\leq
\lambda_1$ for $j=1,\ldots ,N$.
\end{proof}

The error bound of Theorem \ref{as2:thm} now follows from the
estimates
\begin{eqnarray*}
\| P_\lambda -P_\lambda^\infty\|_\Delta &\stackrel{(i)}{\leq} & \|
P_\lambda -P_\lambda^{(m)}\|_\Delta +\|
\mathcal{E}_\lambda^{(m)}\|_\Delta
\\
&\stackrel{(ii)}{\leq} &  |\mathcal{N}_\lambda -1|+ \dim
(\mathcal{A}_\lambda^W) \| \mathcal{E}_\lambda^{(m)}\|_\Delta ,
\end{eqnarray*}
whence $\| P_\lambda -P_\lambda^\infty\|_\Delta= O( \lambda_1^N\,
e^{-\epsilon\, m (\lambda)})$ by Proposition \ref{as-bound:prp},
Proposition \ref{normbound2:prp}, and Lemma \ref{dimbound:lem}.
Here step {\em (i)} hinges on Proposition \ref{as-bound:prp}, and
in step {\em (ii)} we employed Eq. \eqref{Pm-orthogonalexp}
combined with Lemma \ref{termbound:lem}.

\subsection{Proof of Theorem \ref{exact:thm}}
If the reduced $c$-functions $\hat{c}_0 (z)$ and $\hat{c}_1(z)$
are polynomial in $z$ of degree at most $M$ and $2M$,
respectively, then
$P_\lambda^\infty(\mathbf{x})=P_\lambda^{(m)}(\mathbf{x})$ (i.e.
$\mathcal{E}_\lambda^{(m)}=0$) for $m\geq M$. Hence, by
Proposition \ref{triangular:prp}, one has in this case that
\begin{equation}\label{polexp}
P_\lambda^\infty(\mathbf{x}) =
\sum_{\mu\in\Lambda,\,\mu\preceq\lambda} a_{\lambda\mu}^\infty\,
m_\mu (\mathbf{x})\quad ( a_{\lambda\mu}^\infty \in\mathbb{C}),
\end{equation}
provided $m(\lambda)\geq M-1$. Invoking of Proposition
\ref{orthogonality:prp}, and comparing with the defining relations
for $P_\lambda (\mathbf{x})$ in Eqs. \eqref{op1}, \eqref{op2},
leads to the conclusion that the asymptotic functions in question
coincide with the latter polynomials up to normalization
\begin{equation}
P_\lambda (\mathbf{x})= \frac{P_\lambda^\infty(\mathbf{x})}{\|
P_\lambda^\infty \|_\Delta} \qquad \text{for}\;\; m(\lambda) \geq
M-1.
\end{equation}
The (square of the) normalization factor reads
\begin{equation}
\| P_\lambda^\infty \|_\Delta^2
\stackrel{\text{Eq.~\eqref{polexp}}}{=} \langle
P_\lambda^\infty,{\textstyle\sum_{\mu\in\Lambda,\,\mu\preceq\lambda}}\,
a_{\lambda\mu}^\infty\, m_\mu\rangle_\Delta
\stackrel{\text{Prop.~\ref{orthogonality:prp}}}{=} a_{\lambda
\lambda}^{\infty},
\end{equation}
which, by Proposition \ref{triangular:prp}, is equal to $1$ when
$m(\lambda)\geq M $.

\section*{Acknowledgments}
Thanks are due to S.N.M. Ruijsenaars for several helpful
discussions.

\bibliographystyle{amsplain}

\end{document}